\documentclass{amsart}
\date{July 3, 2019}
\usepackage[dvips]{graphicx}
\usepackage{amsmath}

\newcommand{\R}{\boldsymbol{R}}

\renewcommand{\L}{{\mathbf \R_{1}^{n+1}}}
\renewcommand{\phi}{\varphi}

\usepackage[usenames]{color}
\newtheorem{theorem}{Theorem}
\newtheorem{lemma}[theorem]{Lemma}
\newtheorem{fact}[theorem]{Fact}

\theoremstyle{definition}
\newtheorem{definition}[theorem]{Definition}

\theoremstyle{remark}

\renewcommand{\L}{{\mathbf \R_{1}^{n+1}}}

\begin{document}

\title[Bernstein-Type theorem]{%
Bernstein-type theorem for zero mean curvature hypersurfaces 
 without time-like points in Lorentz-Minkowski space
}

\author[Akamine]{S.~Akamine}
\author[Honda]{A.~Honda}
\author[Umehara]{M.~Umehara}
\author[Yamada]{K.~Yamada}


\address[Shintaro Akamine]{%
Graduate School of Mathematics, 
Nagoya University, Chikusa-ku, Nagoya 464-8602, Japan}
\email{s-akamine@math.nagoya-u.ac.jp}

\address[Atsufumi Honda]{%
Department of Applied Mathematics, 
Faculty of Engineering, Yokohama National University,
79-5 Tokiwadai, Hodogaya, Yokohama 240-8501, Japan
}
\email{honda-atsufumi-kp@ynu.ac.jp}

\address[Masaaki Umehara]{%
   Department of Mathematical and Computing Sciences,
   Tokyo Institute of Technology,
   Tokyo 152-8552, Japan}
\email{umehara@is.titech.ac.jp}

\address[Kotaro Yamada]{%
   Department of Mathematics,
   Tokyo Institute of Technology,
   Tokyo 152-8551, Japan}
\email{kotaro@math.titech.ac.jp}

\maketitle
\begin{abstract}
Calabi and Cheng-Yau's Bernstein-type theorem asserts that
{\it an entire zero mean curvature
graph in Lorentz-Minkowski $(n+1)$-space
$\L$
which admits only space-like points is a hyperplane}.
Recently, the third and fourth authors proved 
a line theorem for hypersurfaces at their degenerate  light-like points.
Using this, we give an improvement of the
Bernstein-type theorem, and  we show that {\it an entire zero mean curvature
graph in $\R^{n+1}_1$ consisting only of space-like or
light-like points is a  hyperplane}.
This is a generalization of the first, third and fourth authors'
previous result for $n=2$.
\end{abstract}

\section{Introduction}
For the sake of simplicity,
we abbreviate 
\lq {\bf zero mean curvature}\rq\
to \lq {\bf ZMC}\rq.\
We let $\L$ $(n\ge 2)$ be the Lorentz-Minkowski $(n+1)$-space
of signature $(+\cdots+-)$.
Let $f\colon\Omega\to \R$ be
a $C^2$-differentiable function defined on a
domain $\Omega$ in $\R^n$.
Consider the following two associated functions
$$ 
B_f:=1-f_{x_1}^2-\cdots-f_{x_n}^2, \qquad
A_f:=\sum_{i,j=1}^n 
(B_f\delta_{i,j}+f_{x_i}f_{x_j})f_{x_i,x_j},
$$ 
where $f_{x_i}:=\partial f/\partial x_i$, 
$f_{x_i,x_j}:=\partial^2 f/\partial x_i\partial x_j$ and $\delta_{i,j}$
denotes Kronecker's delta.
The graph $G_f:=\{(x,f(x))\in \L \,;\, x\in \Omega\}$
of a function $f\colon\Omega\to \R$ 
is called {\it zero mean curvature graph} (i.e. ZMC-graph) 
if  $A_f$  vanishes identically on $\Omega$
 (cf. \cite[Appendix B]{UY3}).
Also, the graph  $G_f$  is said to be of {\it constant  mean curvature graph} (i.e. CMC-graph) 
if  $A_f^2-kB_f^3$ vanishes identically on $\Omega$
for a certain constant $k\in \R$.
If $\Omega:=\R^n$, $G_f$ is called an {\it entire graph} in $\R^{n+1}_1$.

\medskip
\begin{definition} \rm
Let $f\colon\Omega\to \R$ be a $C^3$-function.
A point where $B_f>0$ $($resp. $B_f<0$, $B_f=0)$ 
is said to be {\em space-like} $($resp. {\it time-like}, 
{\em light-like}$)$ of $G_f$.
The graph $G_f$ consisting only of space-like points
is said to be {\em space-like}.
If $\Omega$ contains  space-like points
and time-like points, then the graph 
$G_f$ is called of {\em mixed type}.
On the other hand, $G_f$ consisting only
of light-like points is said to be {\em light-like}. 
Moreover, a light-like point $x\in \Omega$
is said to be {\em non-degenerate} if
$\nabla B_f\ne 0$ 
at $x$ and is said to be {\em degenerate}
if $\nabla B_f=0$ at $x$, where $\nabla B_f$
 is the gradient vector field of the function $B_f$.
\end{definition}

\medskip
The following theorem
was 
proved by Calabi \cite{C} 
for $n\leq 4$ and by
Cheng-Yau \cite{CY} for $n\ge 5$:

\medskip
\begin{fact}[Bernstein-type theorem for 
space-like ZMC-hypersurfaces]\label{fact:C}
An entire space-like ZMC-graph is a hyperplane in $\L$.
\end{fact}

\medskip
This is an analogue of the classical Bernstein theorem 
for minimal surfaces in the $(n+1)$-dimensional Euclidean 
space for $n\le 7$. 
The assumption that $G_f$ consists only of space-like points
cannot be removed. In fact, 
\begin{equation}\label{eq:T-Graph}
f_0(x_1,x_2,...,x_{n-1},x_n):=x_n+g(x_1),
\end{equation}
gives an entire graph
without space-like points,
where $g(t)$ ($t\in \R$)
is a $C^\infty$-function of one variable.
A systematic construction of entire ZMC-graphs of mixed type 
is given in \cite{Haw}.
The following fact is known.

\medskip
\begin{fact}[The line theorem for ZMC-hypersurfaces]\label{fact:IL}
Let $G_f$ be the ZMC-graph of a $C^4$-differentiable function $f$ defined on a domain $\Omega \subset \R^n$.
If $o\in \Omega$ is a degenerate light-like point,
then there exists a straight line segment $\sigma \, (\subset \R^n)$ 
passing through $o\in \Omega$ such that
$o$ is not an endpoint of $\sigma$ and
$\Omega\cap \sigma$ consists of degenerate 
light-like points of $G_f$, and $\sigma\ni x\mapsto (x,f(x))$ gives a light-like 
line segment in $\L$. 
\end{fact} 

\medskip
The case $n=2$ was proved by Klyachin \cite{Kl}
under $C^3$-differentiability of $f$,
and was generalized for $n\ge 3$ in the third and fourth
authors' work \cite{UY3}, assuming $C^4$-differentiability of $f$.
Fact \ref{fact:IL} 
holds for the ZMC-graph of arbitrary $C^4$-functions $f$ on $\Omega$
satisfying $(D_{\phi}(f):=)A_f-\phi B_f=0$ for a certain  $C^2$-function
$\phi:\Omega\to \R$, and also holds for real analytic CMC-graphs, see \cite{UY3}.
The purpose of this paper is to prove 
the following improvement of the  Bernstein-type theorem:

\medskip
\noindent
{\bf Theorem A.}
\label{thm:main}
{\it An entire  $C^4$-differentiable ZMC-graph in $\L$ 
which does not admit any time-like points is a hyperplane.}

\medskip
Since the line theorem holds for $C^4$-functions $f$
satisfying $D_{\phi}(f)=0$, Theorem A also holds for such $f$.
When $n=2$, Theorem A follows from
a corollary with slightly stronger assertion 
given 
in 
\cite[Theorem A]{AUY},  which states that
{\it an entire smooth ZMC-graph in $\R^3_1$
which is not a plane
admits a non-degenerate
light-like  point if the set of space-like points is
non-empty}. The authors do not know whether
the corresponding assertion still holds even when 
$n\ge 3$, since a line cannot separate a hyperplane for
$n\ge 3$. It should be remarked that a 
similar statement for area-maximizing ZMC-hypersurfaces
in $\L$ was given in \cite[Theorem F]{E}.

On the other hand, it is well-known that (cf. \cite[Appendix]{UY3}),
$$
A_f=B_f \triangle f -\frac 12 \nabla B_f \star \nabla f
$$
holds, where 
$\triangle f:=\sum_{i=1}^n f_{x_i,x_i}$,  
$\nabla B_f:=((B_f)_{x_1},\dots,(B_f)_{x_n})$
and  \lq$\star$\rq\ is the canonical Euclidean inner 
product of $\R^n$.
So each light-like graph is also a ZMC-graph.
Since it is well-known that (cf. \cite{G}, see also \cite[Corollary B]{UY3})
 the line theorem holds
for light-like hypersurfaces under the 
$C^2$-differentiability of $f$, 
we have the following: 

\medskip
\noindent
{\bf Corollary B.}
\label{cor:main}
{\it An entire $C^2$-differentiable light-like graph in $\R^{n+1}_1$ is a light-like hyperplane.}

\medskip
When $n=2$, this assertion has been proved in the appendix of \cite{AUY}.
The line theorem holds not only for ZMC-hypersurfaces  but also
for real analytic CMC-hypersurfaces (cf. \cite[Theorem E]{UY3}). 
Using the same proof for Theorem A, we also obtain the following:

\medskip
\noindent
{\bf Corollary C.}
\label{cor:main}
{\it 
An entire real analytic $CMC$-graph in $\L$
which has no time-like points and does have
a light-like point must be a light-like hyperplane.
}

\medskip
It is known that  a CMC-graph never changes its causal type,
see \cite[Remark 2.2]{HKKUY}. 
In particular, each light-like point on the graph is
degenerate, and lies on a light-like line, by \cite[Theorem E]{UY3}.
Then Corollary~C follows by the same proof for Theorem A. 
In this corollary, the existence of a light-like point
is crucial, since there are entire space-like CMC-graphs which are not 
ZMC (cf. Treibergs \cite{T}).

\section{Proof of Theorem A}\label{sec1}

We first prepare the following. (As  remarked in \cite{E}, it is due to Bartnik):

\medskip
\begin{lemma}
Let $\bar \Omega$ be the closure of a convex domain in $\R^n$ 
which contains an entire line $l$.
Suppose that $f\colon\bar \Omega\to \R$ is a function such that 
each point of $\bar \Omega$
is space-like or light-like.
If $l$ consists only of light-like points, 
then the graph of $f$ lies in a light-like hyperplane.
\end{lemma}

\medskip
If $\bar \Omega$ is not convex, the assertion fails. Such an example
is given in \cite{AUY2}.

\medskip
\begin{proof}
The proof is the same idea as in 
Ecker \cite[Proposition G]{E}
and  Fernandez-Lopez \cite[Lemma 2.1]{FL}:
Since $l$ consists only of light-like points,
$l\ni x\mapsto (x,f(x))$ gives a light-like line $L$
in $\L$.
Since $\bar \Omega$ consists of only  space-like or light-like points,
we have that
$|x-y|\ge |f(x)-f(y)|$ holds for $x,y\in \bar \Omega$,
where $|z|=\sqrt{\sum_{i=1}^n (z_i)^2}$ for $z=(z_1,...,z_n)\in \R^n$.
For each point $P$ on the line $L\, (\subset \L)$, we set
$$
\bar \Lambda_{P}:=\{Q\in \L\,;\,  (Q-P)\cdot (Q-P)\ge 0\},
$$
which is the closure of the exterior of  the light-cone in $\L$ centered at $P$, 
where the dot denotes the Lorentzian inner product.
We set $F(x):=(x,f(x))$ ($x\in \bar \Omega$).
Since $\bar \Omega$ is convex, any line segment $\overline{PF(x)}$ bounded by
$P$ and $F(x)$ lies in $\bar \Lambda_{P}$. So we have
$F(\bar \Omega)=\bigcup_{x\in \bar \Omega}\overline{PF(x)}\subset \bar \Lambda_{P}$
for each $P\in L$.
Thus, we have
$
F(\bar \Omega)\subset \bigcap_{P\in L} \bar \Lambda_{P}.
$
Since $\bigcap_{P\in L} \bar \Lambda_{P}$ is a light-like hyperplane in $\L$,
we obtain the conclusion.
\end{proof}

\medskip
\noindent
({\it Proof of Theorem A})
Let $G_f$ be the entire ZMC-graph of a $C^4$-differentiable function $f\colon\R^n\to \R$ which does not admit any time-like points.
If there are no light-like points on $\R^n$, then
$f$ is $C^\infty$-differentiable, since $f$ is a solution of
the elliptic partial differential equation $A_f=0$.
Thus the assertion follows from Fact 1.
So we may assume that $G_f$ admits at least one
light-like point $x\in \R^n$.
If $x$ is a non-degenerate light-like point, then
the function $B_f$ changes sign, so $G_f$ must admit
time-like points, a contradiction.
Thus, $x$ is a degenerate light-like point.
By Fact 2, the graph of $f$ contains 
a line segment $\sigma \, (\subset \R^n)$ 
passing through $x$ such that
\begin{itemize}
\item[$\bullet$] 
$x$ is not an endpoint of $\sigma$,
\item[$\bullet$] 
$\sigma$ consists of degenerate 
light-like points of $f$, 
and 
\item[$\bullet$] $\sigma\ni x\mapsto (x,f(x))\in \R^{n+1}_1$ 
gives a light-like line segment in $\L$. 
\end{itemize}
We let $l$ be the complete line containing $\sigma$.
If $\sigma\ne l$, then there exists an endpoint $p$
of $\sigma$ on $l$. Then, $p$ itself must be also a
degenerate light-like point, since $p$ is an accumulation point 
of degenerate light-like points.
Then applying Fact 2 again,  there exists a 
light-like line segment $\sigma'$  passing through $p$
such that $p$ is not an endpoint of $\sigma'$.
Since the ambient metric in $\L$ is Lorentzian,
the nullity of the induced metric of $G_f$ is at most dimension one,
and so the null direction at $p$ is uniquely determined.
In particular, $\sigma'$ also lies in the line $l$.
So the graph $G_f$ contains the entire line $l$
consisting of degenerate light-like points, and
$l\ni x\mapsto (x,f(x))\in \L$ gives a
light-like line.
By Lemma 1,
$f$ must coincide with a light-like hyperplane,
proving the assertion.

\medskip
\noindent
(Acknowledgment)
The authors 
express their gratitude to
Wayne Rossman for helpful comments.

\end{document}